\documentclass[11pt]{article}
\usepackage{latexsym, amsmath, amssymb, amsthm}
\usepackage{enumerate, mathrsfs, bbm}
\usepackage[usenames, dvipsnames]{color}
\usepackage{todonotes}
\usepackage{graphicx, tikz}
\usepackage{authblk, setspace, floatrow}
\usepackage{cite}
\usepackage[top=1 in, bottom=1 in, left=1.2 in, right=1.2 in]{geometry}
\usepackage[pdftex,linkcolor=red,citecolor=blue,backref=page]{hyperref}

\newtheorem{theorem}{Theorem}[section]
\newtheorem{lemma}[theorem]{Lemma}

\newtheorem{remark}[theorem]{Remark}

\numberwithin{equation}{section}

\makeatletter
\def\blfootnote{\xdef\@thefnmark{}\@footnotetext}
\makeatother
\newfloatcommand{capbtabbox}{table}[][\FBwidth]
 
\def\m{\mathbb}						
  \def\eps{\epsilon}	\def\vp{\varphi}
		
\def\p{\partial}  
		
\def\ol{\overline}		
\def\be{\begin{equation}}     \def\ee{\end{equation}}
\def\bp{\begin{pmatrix}}	\def\ep{\end{pmatrix}}

\begin{document}
\title{Smooth solutions to the heat equation which are nowhere analytic in time}
\author[]{Xin Yang,  Chulan Zeng and Qi S. Zhang}
\date{}
\maketitle

\begin{abstract}
The existence of smooth but nowhere analytic functions is well-known (du Bois-Reymond, Math. Ann., 21(1):109-117, 1883).  However, smooth solutions to the heat equation are usually analytic in the space variable. It is also well-known (Kowalevsky, Crelle, 80:1-32, 1875) that a solution to the heat equation may not be time-analytic at $t=0$ even if the initial function is real analytic.  Recently, it was shown in \cite{Zha20, DZ20, DP20} that solutions to the heat equation in the whole space,  or in the half space with zero boundary value,  are analytic in time under essentially optimal conditions. In this paper,  we show that time analyticity is not always true in domains with general boundary conditions or without suitable growth conditions.  More precisely, we construct two bounded solutions to the heat equation in the half plane which are nowhere analytic in time.  In addition,  for any $\delta>0$, we find a solution to the heat equation on the whole plane,  with exponential growth of order $2+\delta$,  which is nowhere analytic in time. 
\end{abstract}




\section{Introduction}
 The study of the existence of nowhere-analytic smooth functions has a rich history (see e.g. \cite{Bil82}) since the pioneering works du Bois-Reymond\cite{dBR1883},  Lerch\cite{Ler1888} and Cellerier\cite{Cel1890}. Later, many other examples were found with different methods,  see e.g. \cite{Bor1895, Had1892, Mor1954, Wal72}.  
For the heat equation, the space analyticity of the classical solution in a space-time domain 
is usually expected as a consequence of parabolic regularity.   But the time analyticity is more delicate and is not true in general, see e.g. the well-known examples in Kowalevsky\cite{Kow1875} and Tychonoff\cite{Tyc1935}.  Under extra assumptions, however,  many time-analyticity results for the heat equation,  Navier-Stokes equations, and some other parabolic equations may still be justified, see e.g. \cite{Wid62,  EMZ17, Gig83, Kom80,Mas67}.   

Recently,  in \cite{Zha20, DZ20}, it was discovered that for any complete and noncompact Riemannnian manifold $M$ whose Ricci curvature is bounded from below,  solutions to the heat equation on $M$ with exponential growth of order 2 are analytic in time.  In particular,  as a corollary to Theorem 2.1 in \cite{DZ20}, for any time interval $(a,b]\subseteq\m{R}$, if $u$ is a smooth solution to the heat equation $\p_{t}u-\p_{x}^2 u=0$ on $\m{R}\times (a,b]$ that satisfies for two positive constants $A_1$ and $A_2$,
\be\label{growth cond for analyticity}
 |u(x,t)|e^{-A_{2}x^2}\leq A_{1},  \quad\forall\, (x,t)\in\m{R}\times (a,b],\ee
then $u$ must be time analytic in $t\in(a,b]$.  The growth restraint (\ref{growth cond for analyticity})  is sharp due to the Tychonoff's non-uniqueness example with suitable modifications (e.g.  see Remark 2.3 in \cite{DZ20} for more details).  Later,  similar phenomena were also found in other types of PDEs\cite{Zen21arXiv,  DZZ21arXiv} and in domains with boundary\cite{DP20}.  In particular,  by denoting $\m{R}^{+}=(0,\infty)$,  Theorem 2.1 in \cite{DP20} implies that for any time interval $(a,b]\subseteq\m{R}$, if $v$ is a smooth solution to the heat equation on $\m{R}^{+}\times(a,b]$ with the Dirichlet boundary condition $v(0,t)=0$ and with the growth constraint $ |v(x,t)|e^{-A_{2}x^2}\leq A_{1}$ for any  $(x,t)\in\m{R}^{+}\times (a,b]$,  then $v$ is time analytic in $t\in(a,b]$. 

The study of analyticity of solutions to PDEs has both a long history,  see e.g. the famous Cauchy-Kowalevsky theorem in \cite{Cau1842, Kow1875}, and many applications, such as the time reversibility, the solvability of backward equations and the control theory.  One particular application 
is about control problems involving heat type equations.  For these problems, it is well known that the set of reachable states, though hard to describe exactly, is just a little larger than the set of those that can be reached by the free heat flow.  But until the papers \cite{Zha20, DZ20},  it's even not clear how to characterize the latter in general (see.  e.g.  the comment on page 1 in \cite{LR20}).  For a precise characterization of the reachable states by the free heat flow, see Corollary 2.2 and Remark 2.5 in \cite{DZ20}.  Later, an explicit formula was derived for the control function by representing solutions with power series in time thanks to the time analyticity,  see Theorem 2.1 in \cite{Zha20arXiv}.


In this paper, however, we discover that the time analyticity is hopeless for general boundary conditions or without suitable growth conditions.  More precisely, we construct solutions to the heat equation on the half space-time plane $\{x\geq 0,\,t\in\m{R}\}$ that satisfy the growth condition (\ref{growth cond for analyticity}) but are nowhere analytic in time.  As a byproduct, we also find a solution to the heat equation on the whole space which is nowhere analytic in time and almost satisfies the growth condition (\ref{growth cond for analyticity}).  This example will demonstrate the sharpness of the growth condition (\ref{growth cond for analyticity}) even if the solution is only required to be analytic in time at a single point.

Denote the space-time domain $\Omega_{1}$ as
\be\label{half plane}
\Omega_{1}=\m{R}^{+}\times\m{R}.\ee
We will construct two bounded solutions to the heat equation on $\ol{\Omega}_1$ which are nowhere analytic in time.  Our first example (\ref{trig series}) can be regarded as an extension to the space-time case of du Bois-Reymond \cite{dBR1883}, which itself is based on the Weierstrass function: a continuous but nowhere differentiable trigonometric series on $\m{R}$.  Our second example (\ref{hk-shift}) takes advantage of the heat kernel $\Phi$ on $\m{R}$, defined as in (\ref{fund soln}), and the method of the condensation of singularities \cite{Han1882, Can1882}.  
\be\label{fund soln}
\Phi(x,t)=\left\{\begin{array}{cll}
(4\pi t)^{-\frac12}\exp\big(-\frac{x^2}{4t}\big) & \text{if} &  x\in\m{R}, \, t>0, \\
0 & \text{if} & x\in\m{R},\, t\leq 0.
\end{array}\right.\ee
Although the construction of $u_2$ is direct via the method of condensation of singularities, 
we remark that the method of constructing $u_1$ in (\ref{trig series}) by the Weierstrass type functions may be more flexible to study other evolutionary PDEs such as the Schr{\"o}dinger equation and the wave equation.

\begin{theorem}\label{Thm, nw-anal on half plane}
Define two functions $u_{1},u_{2}:\ol{\Omega}_{1}\to\m{R}$ by 
\begin{eqnarray}
u_1(x,t) &=& \sum_{k=1}^{\infty} e^{-2^{k}}e^{-2^{k}x}\sin\big( 2^{2k+1}t-2^{k}x \big), \label{trig series}\\
u_{2}(x,t) &=& \sum_{k=1}^{\infty}2^{-k}\Phi(x+1, t-r_{k}),  \label{hk-shift}
\end{eqnarray}
where $\{r_{k}\}_{k=1}^{\infty}$ is an enumeration of all the rational numbers.  Then for $i=1, 2$,
$u_{i}\in C^{\infty}(\ol{\Omega}_1)\cap L^{\infty}(\Omega_1)$ and $u_{i}$ satisfies the heat equation on  $\ol{\Omega}_1$.  However,  for any fixed $x_{0}\in[0,\infty)$, the function $u_{i}(x_0,\cdot)$ is nowhere analytic in $t\in\m{R}$.
\end{theorem}


The functions in Theorem \ref{Thm, nw-anal on half plane} are only defined on $\ol{\Omega}_{1}$.   If we want to construct smooth solutions to the heat equation on the whole plane $\Omega_2:=\m{R}\times\m{R}$,  then the solutions have to break the growth constraint (\ref{growth cond for analyticity}).  In addition, it is well-known that this growth constraint is sharp for everywhere time-analyticity.  More precisely,  for any $\delta>0$,  there exists a solution to the heat equation on $\Omega_2$ which grows slower than $e^{A_2 |x|^{2+\delta}}$ but is not time analytic at some point.  Then it is interesting to investigate the following question: 
 {\em Is the growth condition (\ref{growth cond for analyticity}) sharp for somewhere time-analyticity? }
The next result, which is inspired by (\ref{trig series}), gives a positive answer to this question.  


\begin{theorem}\label{Thm, nw-anal on whole plane}
Let $\Omega_2=\m{R}\times\m{R}$ and $\eps\in(0,1)$.  Define $w_{\eps}:\Omega_2\to\m{R}$ by 
\be\label{unbdd trig series}
w_{\eps}(x,t)=\sum_{k=1}^{\infty}e^{-2^{(1+\eps)k}}e^{-2^{k}x}\sin\big(2^{2k+1}t-2^{k}x \big).\ee
Then $w_{\eps}\in C^{\infty}(\Omega_2)$ and $w_{\eps}$ satisfies the heat equation on $\Omega_2$.  However,  for any fixed $x_{0}\in \m{R}$, the function $w_{\eps}(x_0,\cdot)$ is nowhere analytic in $t\in\m{R}$. Meanwhile, there exist positive constants $A_1$ and $A_2$, which only depend on $\eps$, such that 
\be\label{growth for unbdd trig series}
\sup_{x,t\in\m{R}} |w_{\eps}(x,t)|\exp\big(-A_{2}|x|^{1+\frac{1}{\eps}}\big) \leq A_1.\ee
\end{theorem}

\begin{remark}
For any $\eps\in(0,1)$, the function $w_{\eps}$ in the above theorem, when restricted to $\ol{\Omega}_{1}$ where $x\geq 0$,  is also a bounded nowhere time-analytic solution to the heat equation on $\ol{\Omega}_1$.   Furthermore, for any $\delta>0$, by choosing $\eps=\frac{1}{1+\delta}$,  $w_{\eps}$  is bounded by $A_{1}e^{A_2 |x|^{2+\delta}}$ but is nowhere time-analytic on $\Omega_2$.
\end{remark}


\section{Proofs of Theorems \ref{Thm, nw-anal on half plane} and \ref{Thm, nw-anal on whole plane}}
\label{Sec, proofs}

\subsection{Proof of Theorem \ref{Thm, nw-anal on half plane}}
\begin{itemize}
\item We first study $u_1$.  It is straightforward to check that $u_1\in C^{\infty}(\ol{\Omega}_1)\cap L^{\infty}(\Omega_1)$ and $u_1$ satisfies the heat equation on  $\ol{\Omega}_1$.  Next, for any fixed $x_{0}\geq 0$,  we define
\[h(t)=u_1(x_0,t), \quad \forall\, t\in\m{R}.\]
Then it reduces to  prove that $h$ is not analytic at any point $t_{0}\in\m{R}$. By Cauchy-Hadamard theorem, it suffices to show 
\be\label{radius of conv is 0}
\limsup_{n\to\infty} \bigg(\frac{|h^{(n)}(t_0)|}{n!}\bigg)^{\frac1n}=\infty.\ee

For any integer $m\geq 1$,  the $(2m)^{th}$ and the $(2m+1)^{th}$ derivatives of $h$ at $t_0$ can be written as 
\begin{eqnarray*}
h^{(2m)}(t_0) &=& (-1)^{m}\sum_{k=1}^{\infty}e^{-2^k(1+x_0)} 2^{2m(2k+1)}\sin\big(2^{2k+1}t_0-2^{k}x_0\Big),   \\
h^{(2m+1)}(t_0) &=& (-1)^{m}\sum_{k=1}^{\infty}e^{-2^k(1+x_0)} 2^{(2m+1)(2k+1)}\cos\big(2^{2k+1}t_0-2^{k}x_0\big).
\end{eqnarray*}
For any $N\in\m{Z}^{+}$,  there exists a unique $m_{N}\in\m{Z}^{+}$ such that 
\be\label{choice of m_N}
2^{N}(1+x_0)\leq 4m_{N}< 2^{N}(1+x_0)+4.\ee
Define $F_{N}:\m{Z}^{+}\to\m{R}^{+}$ as 
\be\label{weight}
F_{N}(k)= e^{-2^{k}(1+x_0)}2^{2m_{N}(2k+1)}.\ee
Then 
\begin{eqnarray*}
h^{(2m_{N})}(t_0) &=& (-1)^{m_N}\sum_{k=1}^{\infty}F_{N}(k)\sin\Big(2^{2k+1}t_0-2^{k}x_0\Big),   \\
h^{(2m_{N}+1)}(t_0) &=& (-1)^{m_N}\sum_{k=1}^{\infty}2^{2k+1}F_{N}(k)\cos\Big(2^{2k+1}t_0-2^{k}x_0\Big).
\end{eqnarray*}
By the triangle inequality,
\be\label{triangle ineq}
\begin{split}
\big| h^{(2m_N)}(t_0) \big| & \geq F_{N}(N)\big| \sin\big( 2^{2N+1}t_{0}-2^{N}x_0 \big) \big| - \sum_{k\neq N} F_{N}(k), \\
\big| h^{(2m_{N}+1)}(t_0)\big| & \geq 2^{2N+1} F_{N}(N)\big| \cos\big( 2^{2N+1}t_{0}-2^{N}x_0 \big) \big| - \sum_{k\neq N} 2^{2k+1} F_{N}(k).
\end{split}\ee
Since $|\sin(\theta)|+|\cos(\theta)|\geq 1$ for any $\theta\in\m{R}$, adding the two inequalities in (\ref{triangle ineq}) yields
\be\label{ldd for deri}
\big| h^{(2m_N)}(t_0) \big| + \big| h^{(2m_{N}+1)}(t_0)\big| \geq F_{N}(N)-4\Big( \sum_{k\neq N} 2^{2k} F_{N}(k)\Big).\ee

By direct computation,  it follows from (\ref{weight}) that for any $k\geq 1$, 
\be\label{ratio}
\frac{F_{N}(k+1)}{F_{N}(k)}=\frac{2^{4m_{N}}}{e^{2^{k}(1+x_0)}}=\exp\big[4m_{N}\ln 2-2^{k}(1+x_0) \big].\ee
For any fixed $N$, thanks to the choice (\ref{choice of m_N}) of $m_{N}$ and the fact that $\frac12<\ln 2<1$,  $F_{N}(N)$ is the largest term in the sequence $\{F_{N}(k)\}_{k\geq 1}$.  Moreover,  when $N$ is large enough, $F_{N}(N)$ is much larger than the other terms in the sequence $\{F_{N}(k)\}_{k\geq 1}$.  Actually,  it is not difficult to find a positive constant $N_0$, which only depends on $x_0$, such that 
\be\label{dominating term}
\sum_{k\neq N} 2^{2k} F_{N}(k)\leq \frac{1}{100}F_{N}(N), \quad \forall\, N\geq N_0.\ee 

Plugging (\ref{dominating term}) into (\ref{ldd for deri}) leads to 
\be\label{pos ldd for deri}
\big| h^{(2m_N)}(t_0) \big| + \big| h^{(2m_{N}+1)}(t_0)\big| \geq \frac{1}{2}F_{N}(N),  \quad\forall\, N\geq N_0.\ee
By (\ref{weight}), 
\[F_{N}(N)= e^{-2^{N}(1+x_0)}2^{2m_{N}(2N+1)}\geq e^{-2^{N}(1+x_0)}(2^N)^{4m_{N}}.\]
Reorganizing (\ref{choice of m_N}) gives rise to
\be\label{prop of m_N}
\frac{4(m_{N}-1)}{1+x_0}< 2^{N} \leq \frac{4m_{N}}{1+x_0}.\ee
Consequently, 
\[F_{N}(N)\geq  e^{-4m_{N}}\bigg(\frac{4(m_{N}-1)}{1+x_0} \bigg)^{4m_{N}}\geq 2\bigg(\frac{m_{N}-1}{1+x_0} \bigg)^{4m_{N}}.\]
Thus, for any $N\geq N_0$,  it follows from (\ref{pos ldd for deri}) that
\[\big| h^{(2m_N)}(t_0) \big| + \big| h^{(2m_{N}+1)}(t_0)\big| \geq \bigg(\frac{m_{N}-1}{1+x_0} \bigg)^{4m_{N}}.\]
As a result,
\be\label{ldd in final form}
\frac{\big| h^{(2m_N)}(t_0) \big| + \big| h^{(2m_{N}+1)}(t_0)\big|}{(2m_{N}+1)!}\geq \bigg[\frac{(m_{N}-1)^2}{(1+x_0)^{2}(2m_{N}+1)}\bigg]^{2m_{N}}.\ee
Since $m_{N}\to\infty$ as $N\to \infty$, then (\ref{radius of conv is 0})  follows immediately from (\ref{ldd in final form}).

\item Now we consider $u_2$.  Although the expression (\ref{hk-shift}) looks complicated, the conclusion follows directly from an elegant result in \cite{Wal72}. 

\begin{lemma}\label{Lemma, Walczak}(\cite{Wal72})
Let $\vp$ be a bounded $C^{\infty}$ function which is analytic on $\m{R}\setminus\{0\}$ but not analytic at $0$.  Assume there are positive constants $\delta_0$, $A$ and $L$ such that for any $|t|>A$,
\be\label{Walczak cond}
\sup_{n\geq 0}\, \frac{\big|\p_{t}^{n}\vp(t)\big|}{n!}\, \delta_{0}^{n}<L.\ee
Let $\{a_k\}_{k\geq 1}$ be a sequence of non-zero real numbers such that $\sum\limits_{k=1}^{\infty}|a_k|<\infty$.  Let $\{r_{k}\}_{k\geq 1}$ be an enumeration of all the rational numbers.  Define a function $f:\m{R}\to\m{R}$ by
\[f(t)=\sum_{k=1}^{\infty}a_{k}\vp(t-r_{k}).\]
Then $f\in C^{\infty}(\m{R})$ but $f$ is nowhere analytic on $\m{R}$.
\end{lemma}

Next, we will apply Lemma \ref{Lemma, Walczak} to prove the desired result for $u_2$.  First, we recall that the heat kernel $\Phi$ is defined as in  (\ref{fund soln}).  In addition, by noticing $x+1$ is away from 0 for any $x\geq 0$,  we know for any integer $n\geq 1$, there exists some constant $M_n>0$ such that 
\[ |\p_{x}^{n}\Phi(x+1,t)| + |\p_{t}^{n}\Phi(x+1,t)| \leq M_{n}, \quad \forall\ x\geq 0,\ t\in\m{R}. \]
As a result, $u_2\in C^{\infty}(\ol{\Omega}_1)\cap L^{\infty}(\Omega_1)$ and $u_2$ satisfies the heat equation on  $\ol{\Omega}_1$ since
\begin{align*}
 (\p_{t} - \p_{x}^{2}) u_2 & = (\p_{t} - \p_{x}^2) \bigg( \sum_{k=1}^{\infty}2^{-k}\Phi(x+1, t-r_{k})\bigg) \\
& = \sum_{k=1}^{\infty} 2^{-k} (\p_{t} - \p_{x}^2) \Phi(x+1, t-r_{k})=0.
\end{align*}

Then for any fixed $x_0\geq 0$, define 
\[\vp(t) = \Phi(x_0 + 1,  t), \quad \forall\ t\in\m{R}.\]
According to classical estimates on the heat kernel $\Phi$ (see e.g. formula (3.3) in \cite{Kah1969}),  there exists some constant $C>0$ such that for any $n\in\m{N}$,  
\[ |\p_x^n \Phi(x,t)| \leq \frac{C^{n} n^{n/2}}{t^{(n+1)/2}}\ e^{-\frac{x^2}{8t}},\quad  \forall\, x\in\m{R},\ t> 0.\]
Consequently,
\[  |\p_t^n \Phi(x,t)| = |\p_x^{2n} \Phi(x,t)| \leq \frac{2^{n} C^{2n} n^{n}}{t^{n+\frac12}}\, e^{-\frac{x^2}{8t}},\quad  \forall\, x\in\m{R},\ t> 0.\]
In particular, there exists some constant $C_1>0$ such that for any $n\in\m{N}$,
\be\label{est on vp} |\p_{t}^{n} \vp(t)| = |\p_{t}^{n}\Phi(x_0+1, t)| \leq \frac{C_{1}^{n} n^{n}}{t^{n+\frac12}}\, e^{-(x_0+1)^2/(8t)}, \quad\forall\, t>0. \ee
So by choosing $A=1$ and $\delta_0 = \frac{1}{2C_1}$,  it follows from (\ref{est on vp}) that for any $t>A$, 
\[ \frac{|\p_{t}^{n} \vp(t)|}{n!}\,\delta_0^n \leq \frac{C_1^{n} n^{n} }{n!}\,\frac{1}{(2C_1)^n} =\frac{n^n}{n!}\,\frac{1}{2^{n}}. \]
Thanks to the Sterling formula, we conclude that there exists some constant $L>0$ such that 
\be\label{final est on vp} 
\frac{|\p_{t}^{n} \vp(t)|}{n!}\,\delta_0^n < L, \quad \forall\, t>A. \ee
Noticing that $\p_{t}^{n}\vp(t)=0$ for any $t<0$, so (\ref{final est on vp}) is also valid for $|t|>A$.  Therefore, (\ref{Walczak cond}) is justified for $\vp$.  Finally, in Lemma \ref{Lemma, Walczak},  by setting $a_{k}=\frac{1}{2^{k}}$, we conclude that the function 
\[u_{2}(x_0,\cdot)= \sum_{k=1}^{\infty} \frac{1}{2^{k}}\Phi(x_0+1,\cdot-r_{k}) = \sum_{k=1}^{\infty} a_k \vp(\cdot - r_k) \]
 is nowhere analytic in $t\in\m{R}$.

\end{itemize}

\subsection{Proof of Theorem \ref{Thm, nw-anal on whole plane}}

Fix any $\eps\in (0,1)$,  it is readily seen that $w_{\eps}\in C^{\infty}(\ol{\Omega}_2)$ and $w_{\eps}$ satisfies the heat equation on  $\ol{\Omega}_2$.  

Next, for any fixed $x_{0}\in\m{R}$, we will prove that the function $w_{\eps}(x_0,\cdot)$ is nowhere analytic on $\m{R}$.  Define $h_{\eps}(t)=w_{\eps}(x_0,t)$ for $ t\in\m{R}$. Then it reduces to  prove $h_{\eps}$ is not analytic at any point $t_{0}\in\m{R}$.  By Cauchy-Hadamard theorem, it suffices to show 
\be\label{radius of conv is 0, wp}
\limsup_{n\to\infty} \bigg(\frac{|h_{\eps}^{(n)}(t_0)|}{n!}\bigg)^{\frac1n}=\infty.\ee 
The proof of (\ref{radius of conv is 0, wp}) is similar to that for the function $u_1$ in Theorem \ref{Thm, nw-anal on half plane}, so we will  only sketch the process.  For any large $N$ such that 
$2^{\eps N}\geq 2+|x_0|$, there exists a unique $m_{N}\in \m{Z}^{+}$ such that 
\be\label{choice of m_N, wp}
\big(2^{\eps N}+x_0\big)2^{N}\leq 4m_{N}<\big(2^{\eps N}+x_0\big)2^{N}+4.\ee
Define $F_{N}:\m{Z}^{+}\to\m{R}^{+}$ as 
\be\label{weight, wp}
F_{N}(k)= e^{-2^{(1+\eps)k}}e^{-2^{k}x_0}2^{2m_{N}(2k+1)}.\ee
Then similar to \eqref{ldd for deri}, we have
\[\big| h_{\eps}^{(2m_N)}(t_0) \big| + \big| h_{\eps}^{(2m_{N}+1)}(t_0)\big| \geq F_{N}(N)-4\Big( \sum_{k\neq N} 2^{2k} F_{N}(k)\Big).\]
Thanks to the choice (\ref{choice of m_N, wp}) of $m_{N}$, it is not difficult to find a positive constant $N_0$, which only depends on $x_0$ and $\eps$,  such that 
\be\label{dominating term, wp}
\sum_{k\neq N} 2^{2k} F_{N}(k)\leq \frac{1}{100}F_{N}(N), \quad \forall\, N\geq N_0.\ee
Hence, 
\be\label{pos ldd for deri, wp}
\big| h_{\eps}^{(2m_N)}(t_0) \big| + \big| h_{\eps}^{(2m_{N}+1)}(t_0)\big| \geq \frac{1}{2}F_{N}(N),  \quad\forall\, N\geq N_0.\ee
Reorganizing (\ref{choice of m_N, wp}) leads to
\be\label{prop of m_N, wp}
\frac{4(m_{N}-1)}{2^{\eps N}+x_0}< 2^{N} \leq \frac{4m_{N}}{2^{\eps N}+x_0}.\ee

Based on (\ref{pos ldd for deri, wp}) and (\ref{weight, wp}),  for any $N\geq N_0$,
\begin{align*}
\big| h_{\eps}^{(2m_N)}(t_0) \big| + \big| h_{\eps}^{(2m_{N}+1)}(t_0)\big| &\geq \frac{1}{2}e^{-2^{(1+\eps)N}}e^{-2^{N}x_0}2^{2m_{N}(2N+1)}\\
&\geq e^{-2^{N}(2^{\eps N}+x_0)}(2^{N})^{4m_{N}}.
\end{align*}
Then it follows from (\ref{prop of m_N, wp}) that 
\begin{align*}
\big| h_{\eps}^{(2m_N)}(t_0) \big| + \big| h_{\eps}^{(2m_{N}+1)}(t_0)\big| &\geq e^{-4m_{N}}\bigg(\frac{4(m_{N}-1)}{2^{\eps N}+x_0}\bigg)^{4m_{N}} \geq \bigg(\frac{m_{N}-1}{2^{\eps N}+x_0} \bigg)^{4m_{N}}.
\end{align*}
As a result,
\be\label{ldd in final form, wp}
\frac{\big| h_{\eps}^{(2m_N)}(t_0) \big| + \big| h_{\eps}^{(2m_{N}+1)}(t_0)\big|}{(2m_{N}+1)!}\geq \bigg[\frac{(m_{N}-1)^2}{(2^{\eps N}+x_0)^{2}(2m_{N}+1)}\bigg]^{2m_{N}}.\ee
When $N\to\infty$,  it follows from (\ref{choice of m_N, wp}) that $m_{N}\to\infty$ and 
\[\frac{m_{N}}{(2^{\eps N}+x_0)^2}\sim 2^{(1-\eps)N}\to\infty.\]
Then (\ref{radius of conv is 0, wp}) follows from (\ref{ldd in final form, wp}).

Finally, we need to establish the growth constraint (\ref{growth for unbdd trig series}).  Fix $0<\eps<1$.  Then it suffices to find constants $A_{1}$ and $A_{2}$, which only depend on $\eps$, such that for any $x\geq 0$,
\[\sum_{k=1}^{\infty}\exp\big(-2^{(1+\eps)k}+2^{k}x \big) \leq A_{1}\exp\big(A_{2}\,x^{1+\frac{1}{\eps}}\big).\]
Define $g_{k}(x)=\exp[2^{k}(x-2^{\eps k})]$ for any $k\geq 1$.  Then it reduces to justify 
\be\label{upp for comp series}
\sum_{k=1}^{\infty}g_{k}(x) \leq A_{1}\exp\big(A_{2}\,x^{1+\frac{1}{\eps}}\big).\ee
If $0\leq x\leq 100$, then it is readily seen that the above series in (\ref{upp for comp series}) is uniformly bounded by a constant $B_{1}$ which only depends on $\eps$, so it reduces to consider the case when $x>100$.  Define $K\in\m{Z}^{+}$ to be the unique positive integer such that 
\be\label{choice of K}
\frac{1}{\eps}\log_{2}\Big(\frac{x}{2^{1+\eps}-1}\Big)\leq K < 1+ \frac{1}{\eps}\log_{2}\Big(\frac{x}{2^{1+\eps}-1}\Big),\ee
which can be rewritten as 
\be\label{alt exp for K}
\frac{x}{2^{1+\eps}-1}\leq 2^{\eps K} < \frac{2^{\eps}x}{2^{1+\eps}-1}.\ee
In addition, since $x>100$ and $0<\eps<1$, $K\geq 4/\eps>4$. 
By direct computation, 
\be\label{ratio of g_k}
\frac{g_{k+1}(x)}{g_{k}(x)}=\exp\Big( 2^{k}\big[x-2^{\eps k}(2^{1+\eps}-1)\big] \Big).\ee

\begin{itemize}
\item {\bf Case 1: Estimate of $\sum\limits_{k=1}^{K}g_{k}(x)$.}\\
For any $1\leq k\leq K-1$,  it follows from (\ref{ratio of g_k}) and (\ref{alt exp for K}) that $g_{k+1}(x)\geq g_{k}(x)$. As a consequence, 
\be\label{sum 1 bdd by last term}
\sum\limits_{k=1}^{K}g_{k}(x)\leq K g_{K}(x)\leq K\exp(2^{K}x).\ee
Since $x>100$ and $K>4$, then $K\leq \exp(2^{K}x)$.  Moreover,  we can see from (\ref{sum 1 bdd by last term}) and (\ref{alt exp for K}) that 
\be\label{est for sum 1}
\sum\limits_{k=1}^{K}g_{k}(x)\leq \exp\bigg[2 \bigg(\frac{2^{\eps}x}{2^{1+\eps}-1} \bigg)^{\frac{1}{\eps}} x \bigg]=\exp\big(B_{2}\,x^{1+\frac{1}{\eps}} \big),\ee
where $B_{2}$ is a constant which only depends on $\eps$.

\item {\bf Case 2: Estimate of $\sum\limits_{k=K+1}^{\infty}g_{k}(x)$.}\\
For any $k\geq K$, it follows from \eqref{ratio of g_k} and (\ref{alt exp for K})  that $g_{k+1}(x)\leq g_{k}(x)$. In particular,  by choosing $k=K$ and recalling (\ref{est for sum 1}),  we know
\be\label{est on g_(K+1)}
g_{K+1}\leq g_{K}(x)\leq \exp\big(B_{2}\,x^{1+\frac{1}{\eps}} \big).\ee
From (\ref{alt exp for K}), we have $x\leq (2^{1+\eps}-1)2^{\eps K}$.  Plugging this inequality into  \eqref{ratio of g_k} yields 
\[\frac{g_{k+1}(x)}{g_{k}(x)}\leq \exp\big[ -(2^{1+\eps}-1)2^{k}(2^{\eps k}-2^{\eps K}) \big].\]
So for any $k\geq K+1$, 
\begin{eqnarray*}
\frac{g_{k+1}(x)}{g_{k}(x)} \leq \exp\big[ -(2^{1+\eps}-1)2^{k}(2^{\eps (K+1)}-2^{\eps K}) \big] \leq \exp\big[-2^{k}(2^{\eps}-1)  \big].
\end{eqnarray*}
This implies that for any $k\geq K+1$,  
\begin{eqnarray*}
g_{k+1}(x) &=& \bigg(\prod_{i=K+1}^{k}\frac{g_{i+1}(x)}{g_{i}(x)}\bigg) g_{K+1}(x)\\
&\leq &\exp\big[ -(2^{\eps}-1)2^{K+1}(2^{k-K}-1) \big] g_{K+1}(x) \\
&\leq  &\exp\big[ -(2^{\eps}-1)(2^{k-K}-1) \big] g_{K+1}(x)
\end{eqnarray*}
Thus,  by setting $j=k-K$ and adding $j$ from $1$ to $\infty$,
\be\label{est for sum 2}
\sum_{k=K+2}^{\infty} g_{k}(x) \leq g_{K+1}(x) \sum_{j=1}^{\infty}\exp\big[-(2^{\eps}-1)(2^{j}-1) \big] = B_{3}\,g_{K+1}(x),\ee
where $B_{3}$ is a positive constant which only depends on $\eps$. 
\end{itemize}
Combining (\ref{est for sum 1}), (\ref{est on g_(K+1)}) and (\ref{est for sum 2}) together leads to the desired estimate (\ref{upp for comp series}).

\section*{Acknowledgments} 
We wish to thank Professor Hongjie Dong and Xinghong Pan for helpful discussions.  Q.  S.  Zhang wishes to thank the Simons foundation (grant 710364) for its support.

  \medskip

{\small

}

\bigskip
\bigskip

\thanks{(X. Yang) 
Department of Mathematics, University of California, Riverside, CA 92521, USA}

\thanks{Email: xiny@ucr.edu}

\bigskip

\thanks{(C.  Zeng)
Department of Mathematics, University of California, Riverside, CA 92521, USA}

\thanks{Email: czeng011@ucr.edu}

\bigskip

\thanks{(Q.  S.  Zhang) 
Department of Mathematics, University of California, Riverside, CA 92521, USA}

\thanks{Email: qizhang@math.ucr.edu}

\end{document}